\documentclass[11pt]{amsart}
\usepackage[margin=1.02in]{geometry}
\usepackage{amsmath,amssymb,amsthm,mathtools,microtype}
\usepackage[hidelinks]{hyperref}
\hypersetup{pdftitle={A Quarter-Power Bound for Chowla's Cosine Problem},pdfauthor={Abhishek Shankar}}
\allowdisplaybreaks[2]
\newtheorem{theorem}{Theorem}[section]

\newtheorem{lemma}[theorem]{Lemma}

\theoremstyle{remark}

\newcommand{\T}{\mathbb T}
\newcommand{\Z}{\mathbb Z}
\newcommand{\1}{\mathbf 1}
\newcommand{\E}{\mathbb E}
\newcommand{\TV}{\mathrm{TV}}
\newcommand{\dd}{\,\mathrm d}
\title[Chowla's cosine problem]{A Quarter-Power Bound for Chowla's Cosine Problem}
\author{Abhishek Shankar}
\address{OpenProblem.ai}
\email{abhishek.shankar@openproblem.ai}
\thanks{ORCID: 0009-0002-6552-3622.}
\date{November 2026}
\subjclass[2020]{Primary 42A05; Secondary 11B30, 05C50}
\keywords{Chowla's cosine problem, Fourier idempotents, reversible Markov kernels, additive triples}

\begin{document}
\begin{abstract}
For every nonempty finite set $S$ of $N$ positive integers we prove
\[
 -\min_{x\in\T}\sum_{s\in S}\cos(2\pi sx)
 \gg \left(\frac{N}{\log(2N)}\right)^{1/4}.
\]
The proof uses a reversible phase kernel whose one-step and two-step laws
are close in total variation. Comparison of rows and a height cutoff
produce a partition into pieces of bounded mean. A Fourier rounding estimate
then forces many additive triples to cross the partition.
Two choices of phase density yield bounds in terms of negative mass
and an exact conditional-pair defect, respectively.
A logarithmic refinement gives a log-free bound at low triple density.
At fixed positive triple density the log-free bound is a corollary
of known boundary estimates. We give a separate spectral proof
using an integer difference count and a theorem of Yip.
\end{abstract}
\maketitle

\section{Introduction}

For a nonempty finite set $S\subset\mathbb N$, write
\[
 K(S)=-\min_{x\in\T}\sum_{s\in S}\cos(2\pi sx).
\]
Chowla's conjecture predicts $K(S)\gg |S|^{1/2}$ uniformly
\cite{Chowla}. Bourgain and Ruzsa established subpolynomial
bounds \cite{Bourgain,Ruzsa}. Sanders studied a finite-group
version \cite{Sanders}. Jin, Milojevi\'c, Tomon and Zhang obtained the first
polynomial lower bound, with exponent $1/10-o(1)$ \cite{JMTZ}.
Yip obtained a uniform lower bound with exponent $1/6.43$ \cite{Yip}.
Bedert proved a uniform lower bound with exponent
$1/5-o(1)$ \cite{Bedert}. The exponent $1/5$ without a loss was obtained
in \cite{Shankar} by averaging an additive boundary count.
Bu, Khan, Yi, Jha, Chandupatla and Nagpal obtained a log-free
quarter-power bound for $N$-element subsets of $[1,CN]$, with the constant
depending on $C$, and a uniform lower bound involving the $L^1$ norm
of the cosine sum \cite{Bu}. Here we allow arbitrary frequency span.

We work with a nonempty finite $A=-A\subset\Z\setminus\{0\}$.
Write $n=|A|$, $e_j(x)=\exp(2\pi ijx)$, and
$F=\sum_{a\in A}e_a$, and suppose $F\ge -K$.
All integrals on $\T=\mathbb R/\mathbb Z$ use normalized Haar
measure. We use $f*g(x)=\int f(y)g(x-y)\dd y$ and
$\widehat f(j)=\int f(x)e_{-j}(x)\dd x$.
The idempotence identity $F*F=F$ follows from the coefficients $\1_A$.
Since $F+K$ is nonnegative with a unit Fourier coefficient, $K\ge1$.
The row comparison below uses only a density bound polynomial in $n$.

\begin{theorem}[Negative-part bound]\label{thm:main}
Let $W=F^-=\max(-F,0)$, $\ell=\int W$, and
$D=\|W\|_2^2+\|W*W\|_2^2$. Then
\[
 n\ll D\log(2K)\ll K\ell^3\log(2K)\ll K^4\log(2K).
\]
In particular $K(S)\gg (|S|/\log(2|S|))^{1/4}$ for every nonempty
finite $S\subset\mathbb N$.
\end{theorem}

A second estimate uses $P=K+F$. For $K\ge2$, put $f=P/K$,
$c=f*f$, and define
\[
 H(x,z)=\int\frac{f(s-x)f(s-z)}{c(s)}\dd s,\qquad
 \mathcal E=\iint (1-H(x,z))_+ f(x)f(z)\dd x\dd z.
\]
If $X,Y$ are independent with density $f$, the conditional density of
$X$ given $X+Y=s$ is $f(x)f(s-x)/c(s)$. Thus $H$ is the density
relative to $f(x)f(z)\dd x\dd z$ of two independent draws from this
conditional distribution at a common sum. The parameter $\mathcal E$
is its total variation distance from the independent endpoint law.

\begin{theorem}[Conditional-pair bound]\label{thm:posterior}
For $K\ge2$ one has
\[
 n\ll K^3+K^4\mathcal E\log(2K).
\]
In particular, $\mathcal E\le1$ gives a second derivation of the
uniform $K^4\log(2K)$ bound in Theorem~\ref{thm:main}.
\end{theorem}

Both derivations use the same partition and crossing argument. They
differ in the density used and in the estimate of its transition
defect. The following result identifies two regimes where a logarithm
can be avoided. Put
$T_3=|\{(a,b,c)\in A^3:a+b+c=0\}|=\int F^3$.

\begin{theorem}[Triple-count bounds]\label{thm:triple}
One has
\[
 n\ll K^4+D\log\!\left(2+\frac{KT_3}{n^2}\right).
\]
Consequently $T_3\le Cn^2/K$ implies $n\ll_C K^4$.
For every fixed $\delta>0$, the condition $T_3\ge\delta n^2$
also implies $n\ll_\delta K^4$.
\end{theorem}

\section{A partition from a reversible phase kernel}

Let $q$ be even, continuous, nonnegative, and positive almost
everywhere on $\T$. Set
\[
 C=q*q,\qquad b=C-q,\qquad Z=\int qC,\qquad
 \pi(\dd x)=Z^{-1}q(x)C(x)\dd x.
\]
In particular $C>0$ everywhere. The transition kernel
\[
 T_x(\dd y)=\frac{q(y)q(x-y)}{C(x)}\dd y
\]
is reversible with stationary law $\pi$. Put $J_1=\pi(\dd x)T_x(\dd z)$
and $J_2=\pi(\dd x)T_x^2(\dd z)$, and let
$d=\|J_1-J_2\|_{\TV}$, where $\TV$ is half the $L^1$ distance.
Their densities relative to $\pi\otimes\pi$ are
\[
 r(x,z)=\frac{Zq(x-z)}{C(x)C(z)},\qquad
 r_2(x,z)=\int r(y,z)T_x(\dd y).
\]
Write
$\rho=\int\|T_x-T_y\|_{\TV}\dd J_1(x,y)$.

\begin{lemma}[Comparison of rows by density levels]\label{lem:levels}
For each $z$ let a measurable list of at most $M$ nonnegative
levels be given. If $r^*(x,z)$ is a closest level to $r(x,z)$
and $\varepsilon=\int|r-r^*|\dd(\pi\otimes\pi)$, then
\[
 \rho\ll \varepsilon+Md.
\]
In particular, if $r\le B$ for $B\ge1$, then, for each fixed
$0<\tau<1/2$,
$\rho\ll \tau+O_\tau(1+\log B)d$.
\end{lemma}

\begin{proof}
For fixed $x,z$, draw $Y,Y'$ independently from $T_x$ and put
$X=r(Y,z)$, $X'=r(Y',z)$. Let $X^*,(X')^*$ be the rounded values.
For a threshold $t$ define
\[
 h_t(v,w)=\bigl(\min(v-t,t-w)\bigr)_+,\qquad
 \Gamma_t(X)=\E(X-t)_+-(\E X-t)_+.
\]
The expectation of $h_t(X,X')$ is at most
$\min(\E(X-t)_+,\E(t-X)_+)=\Gamma_t(X)$.
If two levels $v^*,w^*$ are distinct, some midpoint $t$ of
consecutive levels lies in the middle half of the interval
between them. Indeed, take the midpoint of the consecutive gap
containing their arithmetic mean. Thus their difference is bounded
by $4(h_t(v^*,w^*)+h_t(w^*,v^*))$. The maps $h_t$ are
$1$-Lipschitz in the maximum norm. Applying the estimate at
the selected threshold, then summing the nonnegative terms
over all midpoints, gives
\[
 \E|X-X'|\ll \E|X-X^*|+\sum_t\Gamma_t(X).
\]
The case of equal rounded values uses the rounding errors alone.

For a fixed threshold index $t=t(z)$, stationarity identifies
the integral of its Jensen gap against $\pi(\dd x)\pi(\dd z)$ with
\[
 \int \bigl[(r-t)_+-(r_2-t)_+\bigr]\dd(\pi\otimes\pi)
 \le \int(r-r_2)_+\dd(\pi\otimes\pi)=d.
\]
If the number of levels varies, pad the lists by repeats.
Integrating the rounding error costs $\varepsilon$.
By Jensen, the integrated bound for $\E|X-X'|$ also
bounds the integrated $\E|X-\E X|$, which is twice the
average distance between $T_y$ and $T_x^2$ along $J_1$.
The average distance between $T_x^2$ and $T_x$ is $d$.
The triangle inequality proves the first claim.
For the second, use zero and a geometric grid from $\tau$ to
at least $B$ with successive ratio $1+\tau$. The integrated
error is at most $2\tau$, since $\int r=1$.
\end{proof}

\begin{lemma}[Partition into bounded-mean pieces]\label{lem:partition}
There is a finite measurable partition $(E_i)$ of $\T$ with
$k_i=\int_{E_i}q\le2$ and
\[
 \Delta_q=Z-\sum_i\int q_i(q_i*q_i)
 \ll Z\rho+\int q|b|,\qquad q_i=q\1_{E_i}.
\]
\end{lemma}

\begin{proof}
Put $s(x,y)=\|T_x-T_y\|_{\TV}$ and fix $u=1/20$. A point is
good if $T_x(\{y:s(x,y)\le u\})\ge1-u$. By Markov's inequality the bad set has
$\pi$ mass $O(\rho)$. If two good points have distance
greater than $2u$, their $u$-balls are disjoint, and their rows
then have distance at least $1-2u$. Since $4u<1-2u$,
distance at most $2u$ is an equivalence relation on the good
set. Its classes are measurable. There are finitely many because
$x\mapsto T_x$ is continuous in $L^1$ of densities on a compact
domain and distinct classes are separated by $1-2u$.
The $J_1$ mass outside pairs in a common good class is $O(\rho)$.

For a good class $E$ of mass $p>0$, put $\pi_E=\pi|_E/p$,
$\alpha=(\pi_E T)(E^c)$ and
$e=\|\pi_E(\dd x)T_x(\dd y)-\pi_E\otimes\pi_E\|_{\TV}$.
By reversibility the restriction of $\pi_E T$ to $E$ has
density $T_y(E)$ relative to $\pi_E$, so
$\|\pi_E T-\pi_E\|_{\TV}=\alpha\le e$. Averaging row distances
and using the triangle inequality yields
\[
 e\le \int s(x,y)\dd(\pi_E\otimes\pi_E)+\alpha.
\]
The integral is at most $\int_{E^2}s(x,y)\pi_E(\dd x)T_x(\dd y)+2u e$,
since $0\le s\1_{E^2}\le 2u$. Summing over classes with weights
$p$ gives $\sum_E p e\ll\rho$, since
$\sum_E p\alpha=J_1(\bigcup_E E\times E^c)=O(\rho)$ by the
cross-class mass bound, independently of $\alpha\le e$.

Write $Z_E=Zp$, $h(x)=\log(C(x)/Z_E)$ and
$\phi(t)=\int e^{ith}\dd\pi_E$. For independent
$X,Y\sim\pi_E$, the law of $V=X-Y$ is within $2e$ of
$\pi_E$ in total variation. To see this, under
$\pi_E(\dd x)T_x(\dd y)$ the difference has law $\pi_E T$,
by the invariance of the conditional density under $y\mapsto x-y$.
Also
\[
 h(V)-h(X)-h(Y)=\log r_E(X,Y)+\log(C(V)/q(V)),\qquad
 r_E(x,y)=\frac{Z_E q(x-y)}{C(x)C(y)}.
\]
Here $\E|r_E(X,Y)-1|\le2e$, by restriction of the variation
comparison to $E^2$. The logarithms are defined almost surely.
For $|t|\le1$ and $a>0$,
$|e^{it\log a}-1|\ll\min(1,|a-1|)$. Apply this
to $r_E$ and to $q/C$, using the law of $V$, to obtain
\[
 |\phi(t)^2-\phi(t)|\ll
 \delta_E:=e+\int \min(1,|b|/C)\dd\pi_E .
\]
Indeed, the left side compares $\E e^{ith(V)}$ with
$\E e^{it(h(X)+h(Y))}=\phi(t)^2$, up to the $O(e)$ cost of
replacing the law of $V$ by $\pi_E$. The first logarithm costs
$O(e)$ by the bound on $r_E-1$. The second costs
$O(\int\min(1,|b|/C)\dd\pi_E+e)$ by the same law comparison,
using the inverse phase for $C/q$.
For small $\delta_E$, this places $\phi(t)$ in two disjoint
disks of radius $O(\delta_E)$ about $0$ and $1$. Continuity
and $\phi(0)=1$ select the disk about $1$ on $[-1,1]$.
Integration in $t$ then gives
$\int(1-\sin h/h)\dd\pi_E\ll\delta_E$, where the integrand
has its continuous value at zero. It is bounded below by a
positive constant on $|h|\ge\log 2$. Hence
$\pi_E\{C<Z_E/2\}\ll\delta_E$, also for large $\delta_E$
by increasing the constant.

Retain only $E\cap\{C\ge Z_E/2\}$ from each good class. It has
$q$ mass at most two, by $\int_E qC=Z_E$. The discarded
$\pi$ mass is $O(\rho+Z^{-1}\int q|b|)$, counting the bad
set as well. Subdivide the discarded set (the complement of the retained pieces) into finitely many pieces
of $q$ mass at most two using nonatomicity. The probability
that $(X,Y)$ crosses the resulting partition under $J_1$ is
$O(\rho+Z^{-1}\int q|b|)$. Under the triple law with
density $Z^{-1}q(x)q(y)q(x-y)$ the same is true of
$(X,X-Y)$. If the three phases do not all occupy the same piece,
at least one of these two pairs crosses. Multiplying by $Z$
proves the assertion.
\end{proof}

\section{Fourier rounding and crossing mass}

We use the following bounded-mean version of a robust boundary
estimate. Its proof uses the coefficient gap in Bedert's boundary
test and the exact total-boundary inequality from \cite{Shankar}.

\begin{lemma}[Bounded-floor approximation]\label{lem:robust}
Let $G$ be real, even, of mean zero, with $G\ge-2$ almost
everywhere and $|\widehat G(j)|\le2$. For any finite symmetric
$S\subset\Z\setminus\{0\}$,
\[
 |S|\ll 1+\|G-F_S\|_2^2,\qquad F_S=\sum_{a\in S}e_a,
\]
provided $G\in L^2(\T)$.
\end{lemma}

\begin{proof}
We first treat trigonometric polynomials. Let $D_0=\|G-F_S\|_2^2$.
Since $v-v^2$ vanishes at $0,1$ and is Lipschitz on $[-2,2]$,
$\|G-G*G\|_2\ll\sqrt{D_0}$. Write $U_0=G^+$ and $V_0=|G|$.
By expanding positive and negative parts,
\[
 V_0*V_0=G*G+4G^+*G^-\le U_0+|G-G*G|+16.
\]
For $t\ne0$ let
$X=G*(e_tG)$ and $Y=X*(e_{-t}G)$.
Radial clipping at $U_0$ replaces $X$ by $X'$ with
$\|X-X'\|_2\ll1+\sqrt{D_0}$, since $|X|\le V_0*V_0$.
Let $Y'$ be the radial clipping of $X'*(e_{-t}G)$ at $U_0$.
Since $|X'*(e_{-t}G)|\le V_0*V_0$, this costs
$O(1+\sqrt{D_0})$ in $L^2$. By Parseval and the bound of two
on the Fourier multiplier,
\[
 \|(X-X')*(e_{-t}G)\|_2\le 2\|X-X'\|_2.
\]
Thus $\|Y-Y'\|_2\ll 1+\sqrt{D_0}$.
Put $R=X'-Y'$. Its modulus is even, $|R|\le2U_0$, and
$\|R\|_1\le4$. The coefficients of $X-Y$ are
$\widehat G(a)\widehat G(a-t)(1-\widehat G(a+t))$.
Write $d_a=\widehat G(a)-1_S(a)$. Replacing the three
factors one at a time bounds the coefficient error by
$C(|d_a|+|d_{a-t}|+|d_{a+t}|)$, since each unchanged factor
has absolute value at most three. Translation preserves the
$\ell^2$ norm, and $\|d\|_{\ell^2}=\sqrt{D_0}$. Parseval
and the two clipping estimates give
\[
 \|R-F_{B_t}\|_2\ll1+\sqrt{D_0},\qquad
 B_t=\{a\in S:a-t\in S,\ a+t\notin S\}.
\]

For the coefficient test use
$g_t(x)=2(1+\sin(2\pi tx))G(x)$ and
$u_t(x)=(1-\cos(2\pi tx+\pi/4))(g_t*g_t)(x)$.
Let $\psi_t,\varphi_t$ be its even and odd parts.
The mean of $g_t$ is zero and $g_t\ge-8$, so $g_t*g_t\ge-64$
by convolving $g_t+8$ with itself. Hence
$|\varphi_t|\le\psi_t+128$.
For the corresponding test built from $F_S$,
Bedert's Lemma~7.2 \cite{Bedert} gives
$\widehat\psi_t(j)\le c_1=4+1/\sqrt2$ at every frequency and
$\left|\int\varphi_t F_{B_t}\right|\ge c_0|B_t|$ with
$c_0=4+\sqrt2$. The last bound uses the common sign of
the imaginary coefficients on $B_t$. The test coefficients
are uniformly bounded for $G$ as well. They are fixed
quadratic polynomials of $\widehat G(j+at)$, $|a|\le2$,
so the coefficient error from replacing $G$ by $F_S$
is $O(\sqrt{D_0})$ in $\ell^2$. Indeed,
$|uv-u'v'|\le |u||v-v'|+|v'||u-u'|$ bounds each monomial
error by a constant times the sum of the five translated
$|d_j|$. Taking even and odd parts also uses reflected
indices, which preserve the same norm.

Put $L=|R|$. Using that $L$ is even, Parseval and
$\|\widehat L^{\,2}\|_{\ell^2}\le\|L\|_1\|L\|_2$ give
\[
 \int(\psi_t+128)(L*L)
 \le c_1\|R\|_2^2+O(1+\sqrt{D_0}\|R\|_2).
\]
The same integral is at least $|\int\varphi_t(R*R)|$.
The uniform coefficient bound controls the replacement of
$R*R$ by $F_{B_t}$ at cost
$O(\sqrt{|B_t|}\|R-F_{B_t}\|_2+\|R-F_{B_t}\|_2^2)$.
The error of the test coefficients on $B_t$ costs
$O(\sqrt{D_0|B_t|})$. Comparing upper and lower estimates
and writing $b_t=|B_t|$ gives
\[
 (c_0-c_1)b_t\ll (1+\sqrt{D_0})\sqrt{b_t}+1+D_0,
\]
since $\|R\|_2\le\sqrt{b_t}+O(1+\sqrt{D_0})$.
Young's inequality absorbs the square-root term.

Let $m=|S|$. We may assume $m$ exceeds a sufficiently large
absolute multiple of $1+D_0$. Cauchy--Schwarz for the
measure $(G+2)\dd x$ gives
$\int (G+2)F_S^2\ge (m-\sqrt{mD_0})^2/2$.
Since $\|F_S^2\|_2\le m^{3/2}$, it follows that
$\int F_S^3\gg m^2$. The total-boundary inequality
$\int F_S^3\le3\sum_{t\in S}|B_t|$ \cite{Shankar}
now yields $m\ll1+D_0$. For general $G\in L^2$, convolve with
mean-one Fej\'er kernels and take the $L^2$ limit with $S$
fixed. Positivity of the kernels preserves the floor and
the coefficient bound.
\end{proof}

\begin{lemma}[Crossing inequality]\label{lem:crossing}
For any partition of $q$ as in Lemma~\ref{lem:partition} with
$k_i\le2$, writing $C_i=q_i*q_i$, one has
\[
 \|q\|_2^2-\left(\int q\right)^2
 \ll \int q+\Delta_q+\int q b_-+\sum_i\int C_i b_+.
\]
\end{lemma}

\begin{proof}
We have $\sum_i C_i\le C$. On $E_i$ the negative squared error
of $C_i-q_i$ is at most $q(C-C_i)+qb_-$, and the positive
squared error is at most $C_i b_+$. Off $E_i$ we have
$C_i^2\le C_i(q+b_+)$. Summing gives
\[
 \sum_i\|C_i-q_i\|_2^2\le
 2\Delta_q+\int qb_-+\sum_i\int C_i b_+.
\]
Here $\sum_i\int_{E_i^c}qC_i\le\Delta_q$ follows from
$\sum_i C_i\le C$, and $\sum_i\int_{E_i}q(C-C_i)=\Delta_q$.
For each $i$, round $\widehat q_i(j)$ at $j\ne0$ to the
nearest of $0$ and $1$, with the same tie rule at $\pm j$.
The set $S_i$ rounded to one is finite by square summability
and symmetric by conjugation symmetry. For any complex $z$
the distance to $\{0,1\}$ is at most $2|z^2-z|$.
Thus the squared error of rounding is at most
$4\|C_i-q_i\|_2^2$ by Parseval.
The even part of $q_i-k_i$ has floor at least $-2$ and
coefficients bounded by two. Taking the even part can only
decrease its squared distance to $F_{S_i}$. Apply
Lemma~\ref{lem:robust} to each nonempty $S_i$.
There are at most $2\sum_i k_i$ such sets because each needs
$k_i\ge1/2$. Finally
$\sum_i\sum_{j\ne0}|\widehat q_i(j)|^2=\|q\|_2^2-\sum_i k_i^2
\ge\|q\|_2^2-(\int q)^2$. Bounding this energy by twice
the sum of squared rounding errors plus $2\sum_i|S_i|$
gives the result.
\end{proof}

\section{Two phase densities}

Combining Lemmas~\ref{lem:partition} and~\ref{lem:crossing}
bounds the nonconstant Fourier energy of $q$ by
\begin{equation}\label{eq:master}
 \|q\|_2^2-\left(\int q\right)^2
 \ll \int q+Z\rho+\int q|b|+\sum_i\int C_i b_+.
\end{equation}
Here $(E_i)$ is the partition produced by
Lemma~\ref{lem:partition}. We retain the sum of $C_i$ rather
than $C$ in the final term because it will give a better
bound for $q=P$.

\subsection{The positive part of the Fourier sum}

Put $U=F^+$ and take $q=Q=U+\ell$, with $\ell,D$
as in Theorem~\ref{thm:main}. Since $\int U=\int W=\ell$,
a unit Fourier coefficient gives $2\ell=\|F\|_1\ge1$.
We also have $\ell\le K$ and $C=U*U+3\ell^2\ge3\ell^2$.
The identity $F*F=F$ gives
\begin{equation}\label{eq:b}
 b=2U*W-W*W-W+3\ell^2-\ell.
\end{equation}
We claim
\begin{equation}\label{eq:d-bound}
 \|b\|_2^2+\int C|b|+\int Q|b|\ll D,
 \qquad Z=n+O(D).
\end{equation}
Let $u_j=\widehat U(j)$ and $w_j=\widehat W(j)$.
For $j\ne0$, $u_j=\1_A(j)+w_j$ and
$\widehat b(j)=(2\1_A(j)-1)w_j+w_j^2$.
The zero coefficient is $4\ell^2-2\ell$.
Thus $\|b\|_2^2\ll D$ since $\ell^4\le\|W*W\|_2^2$.
To bound $\int C|b|$, majorize the absolute value in
\eqref{eq:b} term by term. In the three convolution pairings,
Parseval gives
\[
 \begin{aligned}
 \int C(U*W)&=\sum_j u_j^3w_j+3\ell^4,\\
 \int C(W*W)&=\sum_j u_j^2 w_j^2+3\ell^4,\\
 \int CW&=\sum_j u_j^2 w_j+3\ell^3.
 \end{aligned}
\]
The coefficients are real by evenness. In the expansions
the linear sum is $\sum_{j\in A} w_j=\int FW=-\|W\|_2^2$.
All terms of degrees two through four are bounded by a
constant times $\sum_j(|w_j|^2+|w_j|^4)=D$.
The remaining constants are $O(D)$ because
$\int C=4\ell^2$ and $\ell\ge1/2$.
It follows that $\int C|b|\ll D$ and
$\int Q|b|\le\int C|b|+\|b\|_2^2\ll D$.
Finally $\int Q^2=n-\|W\|_2^2+3\ell^2$ and
$|Z-\int Q^2|\le\int Q|b|$, proving \eqref{eq:d-bound}.

For the transition defect, the density of $J_2-J_1$ with
respect to Haar measure on $\T^2$ is
\begin{equation}\label{eq:two-step}
 \frac{q(x)q(z)}{Z}
 \left[b(x-z)-\int q(x-y)q(z-y)\frac{b(y)}{C(y)}\dd y\right].
\end{equation}
Indeed $q/C=1-b/C$ in the integral formula for the
two-step density, and the term without $b/C$ integrates
to $C(x-z)$. In \eqref{eq:two-step} keep $b_+$ in the first
term and $b_-$ in the second to bound the positive part.
The integral of the resulting majorant equals
$Z^{-1}\int C|b|$. Hence $d\ll D/Z$ for $Q$.

The bounds $\|Q\|_\infty\le n+\ell$ and $C\ge3\ell^2$ show that the one-step density $r$ relative
to $\pi\otimes\pi$ is bounded by a fixed polynomial in $n$.
For example $Z\le(\int Q)^2\|Q\|_\infty$.
Lemma~\ref{lem:levels} gives
$\rho\ll\tau+O_\tau(D\log(2n)/Z)$ for any fixed small $\tau>0$.
The nonconstant Fourier energy is $n-O(D)$ and $\int Q\ll D$,
since $\ell^4\le D$ and $\ell\ge1/2$.
Using \eqref{eq:d-bound} in \eqref{eq:master} gives
$n\ll D+Z\rho$. Choose $\tau$ as a sufficiently small
absolute constant and use $Z=n+O(D)$ to absorb its contribution.
This first gives $n\ll D\log(2n)$.
Young's inequality and $W\le K$ give
\[
 D\le(1+\ell^2)\|W\|_2^2\le(1+\ell^2)K\ell\ll K\ell^3.
\]
In particular $n/\log(2n)\ll K^4$. Since $\log(2n)\ll\sqrt n$,
this implies $n\ll K^8$ and hence $\log(2n)\ll\log(2K)$.
Substituting into the negative-part estimate gives its stated form.
For the cosine formulation take $A=S\cup(-S)$ and
$K=2K(S)\le2|S|$. This completes the proof of
Theorem~\ref{thm:main}.

\subsection{The exact conditional-pair defect}

Take $q=P=K+F$, with $K\ge2$. The polynomial is positive
almost everywhere, and $C=P+a_*$ where $a_*=K^2-K$.
Also $Z=n+K^3$, $c=C/K^2$, and $b=a_*$ is constant.
Equation~\eqref{eq:two-step} becomes the signed measure
\[
 J_2-J_1=\frac{a_*K^2}{Z}(1-H(x,z)) f(x)f(z)\dd x\dd z .
\]
In particular the total variation distance between the two laws is exactly
$d=a_*K^2\mathcal E/Z$. The fact that $H$ is a probability
density relative to the independent endpoint law follows
also from $\int f(x)f(s-x)\dd x=c(s)$.

The one-step density is again bounded by a polynomial in $K$,
using $C\ge a_*$ and the bound $n\ll K^8$ just proved. Lemma~\ref{lem:levels}
therefore gives $\rho\ll\tau+O_\tau(d\log(2K))$.
For the partition with $k_i\le2$ we have
\[
 \int P|b|=a_*K,\qquad
 \sum_i\int C_i b_+=a_*\sum_i k_i^2\le 2a_*K .
\]
The nonconstant energy of $P$ is $n$. It follows from
\eqref{eq:master} that $n\ll K^3+Z\rho$.
Absorbing the fixed small $\tau$ term as above and using
$Zd=a_*K^2\mathcal E\le K^4\mathcal E$ proves
Theorem~\ref{thm:posterior}.

\section{A refinement by the number of additive triples}

We return to $q=Q=F^++\ell$ and its associated $C,Z,J_1$.
Write $B=\int Q^2$, let $(x,z)\sim J_1$, and set
$s=x-z$, $R=Q(s)/Q(x)$. Under $J_1$ both $x$ and $s$
have law $\pi$. The joint density in $x,s$ is
$Q(x)Q(s)Q(x-s)/Z$, so $R$ and $R^{-1}$ have the same
law. Let $\Pi_A$ be orthogonal Fourier projection onto
$A$. We have
\begin{equation}\label{eq:moment}
 Z\E_{J_1}\sqrt R
 =\langle \Pi_A Q^{3/2},\Pi_A Q^{1/2}\rangle+O(K^2\sqrt B).
\end{equation}
Indeed the numerator is
$\iint Q(s)^{3/2}Q(x)^{1/2}Q(x-s)\dd x\dd s$.
Replacing its last factor by $F(x-s)$ gives the projection
pairing. Its cost is at most
$\|Q-F\|_\infty(\int Q^{3/2})(\int Q^{1/2})\ll K^2\sqrt B$
by Cauchy--Schwarz, since $Q-F=W+\ell$ and $\int Q\le2K$.

It suffices to consider $n$ exceeding a sufficiently large
absolute multiple of $K^4$. In this range $B,Z\asymp n$
by \eqref{eq:d-bound} and $D\ll K^4$. Also
$\int Q^3=T_3+O(Kn)$, by factoring the difference of cubes
and using $\|Q-F\|_\infty\le2K$ and
$\int(Q^2+F^2)\ll n$. Applying Cauchy--Schwarz in
\eqref{eq:moment} gives
\[
 \E_{J_1}\sqrt R\ll 1+\frac{\sqrt{KT_3}}n.
\]
Let $L=A_0(1+\sqrt{KT_3}/n)^2$ with a sufficiently large
absolute $A_0$. By Markov's inequality and the reciprocal
symmetry, $R\in[L^{-1},L]$ outside an arbitrarily small
fixed $J_1$ mass, with the choice of $A_0$ depending only
on that mass. Also $C(x)/Q(x)\in[1/2,2]$ outside
$O(D/Z)$ stationary mass: on its complement
$C\le2|C-Q|$, and \eqref{eq:d-bound} applies.

Thus $Q(x-z)/C(x)$ lies in $[(2L)^{-1},2L]$ outside
small $J_1$ mass. Use a geometric grid on that interval
with successive ratio at most $1+\tau$, multiply the levels
by $Z/C(z)$, and include zero. These are the levels of
Lemma~\ref{lem:levels} for $r(x,z)$. The integrated error is
at most $\tau$ plus the $J_1$ mass of the exceptional pairs,
because zero is an available level. The number of levels is
$O_\tau(1+\log L)$. Accordingly
\[
 \rho\ll \tau+J_1\{R\notin[L^{-1},L]\}
               +O_\tau\bigl((1+\log L)D/Z\bigr).
\]
Use $n\ll D+Z\rho$, already shown for $Q$, and absorb
the first two terms by choosing $\tau$ and the exceptional
mass sufficiently small. Since
$1+\log L\ll\log(2+KT_3/n^2)$ for the fixed $A_0$,
this yields $n\ll D\log(2+KT_3/n^2)$ in the large range.
Including $O(K^4)$ handles the other range and proves
the first part of Theorem~\ref{thm:triple}.

\section{The graph argument at fixed positive triple density}

Let $\mathsf A$ be the adjacency matrix of the graph on $A$
where $a,b$ are joined if $a-b\in A$. There are $T_3$
ordered edges. The matrix $\mathsf A+KI$ is positive
semidefinite because for any vector $v$ indexed by $A$ its
quadratic form is $\int P(x)|\sum_{a\in A}v_a e_a(x)|^2\dd x$.

For a proposed clique block $V\subset A$ of size $m$, let
$p,q$ be the numbers of its positive and negative elements.
The number of ordered pairs in $V^2$ whose difference lies
in $V$ equals the number of ordered pairs in $V^2$ whose sum
lies in $V$. For each positive summand the possible positive
targets are larger, so positive-positive pairs contribute at
most $p(p-1)/2$. Ordering the absolute values of the
negative elements gives $q(q-1)/2$ for negative-negative
pairs. Opposite signs contribute at most $2pq$.
The total is at most $3m^2/4-m/2$, since $pq\le m^2/4$.
Consequently for $m\ge3$ at least $m(m-1)/8$ distinct
ordered pairs $(a,b)$ have $a-b\notin V$.
The same is true for $m=2$: a set of two nonzero integers
cannot contain both nonzero differences of its elements.

For a partition into proposed clique blocks, each such bad
pair is either a missing edge within a block or, when
$a-b\in A\setminus V$, gives the existing cross-block edge
$(a,a-b)$. This map on the second type is injective and the
two types give distinct ordered edits. If $L$ is the number
of proposed ordered clique edges and $E$ is the ordered edit
cost, then $L\le8E$ and $T_3\le L+E\le9E$.
In particular the unordered edit distance from every union
of cliques is at least $T_3/18$.

The fixed-density conclusion also follows from earlier boundary
estimates. The total-boundary inequality \cite[Lemma~3.1]{Shankar}
and the bound of Bu et al. \cite[Proposition~9.1]{Bu}, refining
Bedert \cite[Proposition~7.3]{Bedert}, give
\[
 T_3\le 3\sum_{t\in A}|B_t|\ll nK^3\|F\|_1\ll nK^4,
 \qquad B_t=\{a\in A:a-t\in A,\ a+t\notin A\}.
\]
We give a separate spectral proof.

When $T_3\ge\delta n^2$, the graph is $\delta/18$-far from every union of cliques in Yip's convention of at least $\epsilon n^2$ unordered edits. Yip's Theorem~1.7 \cite{Yip}
therefore gives $-\lambda_{\min}(\mathsf A)\gg_\delta n^{1/4}$.
Positivity of $\mathsf A+KI$ bounds this quantity by $K$.
This proves the fixed-density assertion in
Theorem~\ref{thm:triple}, independently of the transition
argument. In particular, on any sequence with $n/K^4\to\infty$,
the triple density $T_3/n^2$ must tend to zero while
$KT_3/n^2$ tends to infinity.


\begin{thebibliography}{99}
\bibitem{Bourgain}
J.~Bourgain, \emph{Sur le minimum d'une somme de cosinus},
Acta Arith. \textbf{45} (1986), 381--389.
\bibitem{Ruzsa}
I.~Z. Ruzsa, \emph{Negative values of cosine sums},
Acta Arith. \textbf{111} (2004), 179--186.
\bibitem{Sanders}
T.~Sanders, \emph{Chowla's cosine problem},
Israel J. Math. \textbf{179} (2010), 1--28.
\bibitem{Bedert}
B.~Bedert, \emph{Polynomial bounds for the Chowla Cosine Problem},
\href{https://arxiv.org/abs/2509.05260v3}{arXiv:2509.05260v3} (2026).

\bibitem{Bu}
B.~Bu, F.~Khan, H.~Yi, A.~Jha, I.~Chandupatla and M.~Nagpal,
\emph{A half-constant upper bound for Chowla's cosine problem},
\href{https://arxiv.org/abs/2610.17160v3}{arXiv:2610.17160v3} (2026).

\bibitem{Chowla}
S.~Chowla, \emph{Some applications of a method of A. Selberg},
J. Reine Angew. Math. \textbf{217} (1965), 128--132.

\bibitem{JMTZ}
Z.~Jin, A.~Milojevi\'c, I.~Tomon and S.~Zhang,
\emph{From small eigenvalues to large cuts, and Chowla's cosine problem},
\href{https://arxiv.org/abs/2509.03490v2}{arXiv:2509.03490v2} (2025).

\bibitem{Shankar}
A.~Shankar, \emph{A log-free $n^{1/5}$ bound for Chowla's cosine problem},
\href{https://arxiv.org/abs/2609.05338}{arXiv:2609.05338} (2026).

\bibitem{Yip}
F.~Yip, \emph{Lower bounds for the magnitude of the minimum
eigenvalue of graphs with applications to MaxCut and Chowla's cosine problem},
\href{https://arxiv.org/abs/2610.01657v2}{arXiv:2610.01657v2} (2026).
\end{thebibliography}
\end{document}